\documentclass[reqno,11pt]{article}
\usepackage{mathrsfs}
\usepackage{amssymb}
\usepackage[usenames,dvipsnames]{xcolor}
\usepackage{amsthm}
\usepackage{amsmath}
\usepackage{amsfonts}
\usepackage{enumerate}

\usepackage{graphicx}
\usepackage{txfonts}

 \theoremstyle{definition}

 \numberwithin{equation}{section}

\newtheorem{theorem}{Theorem}[section]
\newtheorem{lemma}[theorem]{Lemma}
\newtheorem{corollary}[theorem]{Corollary}
\newtheorem{proposition}[theorem]{Proposition}

\theoremstyle{definition}
\newtheorem{definition}[theorem]{Definition}

\theoremstyle{remark}
\newtheorem{remark}[theorem]{Remark}

\begin{document}
\title{On analyticity and temporal decay rates of solutions to the viscous resistive Hall-MHD system}
\author{Shangkun WENG\footnote{Pohang Mathematics Institute,
Pohang University of Science and Technology. Hyoja-Dong San 31, Nam-Gu
Pohang, Gyungbuk 790-784, Korea. {\it Email: skweng@postech.ac.kr, skwengmath@gmail.com}}}
\date{Pohang Mathematics Institute, POSTECH}
\maketitle

\def\be{\begin{eqnarray}}
\def\ee{\end{eqnarray}}
\def\p{\partial}
\def\no{\nonumber}
\def\e{\epsilon}
\def\de{\delta}
\def\De{\Delta}
\def\om{\omega}
\def\Om{\Omega}
\def\f{\frac}
\def\th{\theta}
\def\la{\lambda}
\def\b{\bigg}
\def\al{\alpha}
\def\La{\Lambda}
\def\ga{\gamma}
\def\Ga{\Gamma}
\def\ti{\tilde}
\def\Th{\Theta}
\def\si{\sigma}
\def\Si{\Sigma}
\def\bt{\begin{theorem}}
\def\et{\end{theorem}}
\def\bc{\begin{corollary}}
\def\ec{\end{corollary}}
\def\bl{\begin{lemma}}
\def\el{\end{lemma}}
\def\bp{\begin{proposition}}
\def\ep{\end{proposition}}
\def\br{\begin{remark}}
\def\er{\end{remark}}
\def\bd{\begin{definition}}
\def\ed{\end{definition}}
\def\bpf{\begin{proof}}
\def\epf{\end{proof}}
\begin{abstract}
  We address the analyticity and large time decay rates for strong solutions of the Hall-MHD equations. By Gevrey estimates, we show that the strong solution with small initial date in $H^r(\mathbb{R}^3)$ with $r>\f 52$ becomes analytic immediately after $t>0$, and the radius of analyticity will grow like $\sqrt{t}$ in time. Upper and lower bounds on the decay of higher order derivatives are also obtained, which extends the previous work by Chae and Schonbek (J. Differential Equations 255 (2013), 3971--3982).
\end{abstract}

\begin{center}
\begin{minipage}{5.5in}

\end{minipage}
\end{center}

\section{Introduction and main results}

In this paper we address the analyticity of strong solutions to the incompressible viscous resistive Hall-Magnetohydrodynamic equations.  The incompressible viscous resistive Hall-MHD equations take the following form:
\be\label{hmhd}\begin{cases}
\p_t u+u\cdot\nabla u+\nabla \pi=B\cdot\nabla B +\Delta u,\\
\p_t B-\nabla\times (u\times B)+\nabla\times ((\nabla\times B)\times B)=\Delta B,\\
\text{div} u=\text{div} B=0,
\end{cases}
\ee
where $u(x,t)=(u_1(x,t), u_2(x,t), u_3(x,t))$ and $B(x,t)=(B_1(x,t), B_2(x,t), B_3(x,t))$, $(x,t)\in \Bbb R^3\times[0, \infty)$, are the fluid velocity and magnetic field, $\pi=p+\f{1}{2} |B|^2$, where $p$ is the pressure. We will consider the Cauchy problem for (\ref{hmhd}), so we prescribe the initial data
\be\no
u(x,0)=u_0(x),\quad B(x,0)=B_0(x).
\ee
The initial data $u_0$ and $B_0$  satisfy the divergence free condition,
\be\no
\text{div}\, u_0(x)=\text{div} \,B_0(x)=0.
\ee
The application of Hall-MHD equations is mainly from the understanding of magnetic reconnection phenomena \cite{hg,lighthill,pm}, where the topology structure of the magnetic field changes dramatically and the Hall effect must be included to get a correct description of this physical process. The authors in \cite{adfl} had derived the Hall-MHD equations from a two-fluids Euler-Maxwell system for electrons and ions by some scaling limit arguments. They also provided a kinetic formulation for the Hall-MHD. Recently, there are many researches on the Hall-MHD equations, concerning global weak solutions \cite{adfl,cdl}, local and global (small) strong solutions \cite{cdl,cwu,cl13,fhn}, singularity formation in Hall-MHD \cite{cweng}, and the asymptotic behavior of weak and strong solutions \cite{cs,weng}. In \cite{cweng} we have showed that the Hall-MHD (\ref{hmhd}) without resisitivity is not globally in time well-posed in any $H^m(\mathbb{R}^3)$ with $m>\f 72$, i.e. for some axisymmetric smooth data, either the solution will become singular instantaneously, or the solution blows up in finite time. Note that singularity formation in compressible fluid was proved long time ago, however, whether the incompressible Navier-Stokes equation will develop singularity in finite time is still greatly open. So we believe the Hall-MHD model is not only physical importance but also mathematical interesting, since it provides an example where singularity may delevop in incompressible fluids as shown in \cite{cweng}.

Chae and Schonbek \cite{cs} investigated the temporal decay estimates for weak solutions to Hall-MHD system with initial data in $L^1\cap L^2$. They also obtained algebraic decay rates for higher order Sobolev norms of strong solutions to (\ref{hmhd}) with small initial data. It turned out that the Hall term does not affect the time asymptotic behavior, and the time decay rates behaved like those of the corresponding heat equation. Here we generalized their results to cover more classes of initial data. The proof follows the Fourier splitting method developed in Schonbek and many other authors, one may refer to \cite{bv07b,s1,s2,s4,sss,sw,wiegner} and the reference therein.

Consider the heat system with same initial data $(u_0,B_0)$
\be\label{h}
\begin{array}{lll}
\p_t v= \Delta v,\ \ &v(x,0)=u_0(x),\\
\p_t w= \Delta w,\ \  &w(x,0)= B_0(x).
\end{array}
\ee
Before introducing the result, we define some notations. $\|\cdot\|_p (1\leq p\leq \infty)$ denotes the usual $L^p(\mathbb{R}^3)$ norm. Let $\mathcal{V}=\{v\in (C_0^{\infty}(\mathbb{R}^3))^3: \nabla \cdot v=0\}$ and $H$ be the closure of $\mathcal{V}$ in $(L^2(\mathbb{R}^3))^3$. We also introduce the following weighted function space $W_2=\{v: \|v\|_{W_2}^2:=\int_{\mathbb{R}^3} |x||v(x)|^2 dx<\infty\}$.
\bt\label{m2}(Upper bound).
{\it
Let $(u_0,B_0)\in H\times H$ and $(u(x,t),B(x,t))$ be a weak solution of the Hall-MHD equations with initial datum $(u(x,0),B(x,0))=(u_0(x),B_0(x))$.
\begin{enumerate}[(1)]
  \item Assume that  the solution $(v,w)$ of (\ref{h}) satisfies
  \be\label{m201}
  \|v(t)\|_2^2+ \|w(t)\|_2^2 \leq C(t+1)^{-\alpha}
  \ee
  for all $t\geq 0$, some constant $C>0$ and $\al\geq 0$. Then
  \be\label{m202}
  \|u(t)\|_2^2 +\|B(t)\|_2^2 \leq C(t+1)^{-\bar{\al}},\quad \bar{\al}=\min\{\alpha,\f 52\}.
  \ee
  \item If $0\leq \alpha\leq 5/2$, then there is a constant $C$, depending only on the $L^2$-norm of the initial datum $(u_0,B_0)$ such that for $D(x,t)=(u-v, B-w)(x,t)$, we have
  \be\label{m203}
  \|D(t)\|_2^2 &\leq& \begin{cases}
  \begin{array}{ll}
  C (t+1)^{-5/2},\  &\text{if}\ 1<\alpha\leq 5/2,\\
  C(t+1)^{-5/2}(1+\log^2(t+1)),\ &\text{if}\ \alpha=1,\\
  C(t+1)^{-5/2+2(1-\alpha)},\ &\text{if}\  0\leq \alpha<1.
  \end{array}
  \end{cases}
  \ee
\end{enumerate}
}
\et
Following the ideas developed in \cite{s4,sss}, we also investigate the lower bounds of large time decay rates for weak solutions to the Hall-MHD equations (\ref{hmhd}). Given $u=(u_1,u_2,u_3)$ and $B=(B_1, B_2, B_3)$ in $[L^1(0,\infty; L^2(\mathbb{R}^3))]^3$, introducing the matrices $\widetilde{\mathcal{A}}=[\widetilde{\mathcal{A}}_{ij}]$, $\widetilde{\mathcal{C}}=[\widetilde{\mathcal{C}}_{ij}]$, and $\langle x, B_0\rangle= [\langle x, B_0\rangle_{ij}]$, where
\be\no
\widetilde{\mathcal{A}}_{ij}&=&\int_0^{\infty}\int_{\mathbb{R}^3}(u_i u_j -B_i B_j)(x,t) dx dt,\\\no
\widetilde{\mathcal{C}}_{ij}&=&\int_0^{\infty}\int_{\mathbb{R}^3}(u_i B_j -B_i u_j)(x,t) dx dt,\\\no
\langle x, B_0\rangle_{ij}  &=&\int_{\mathbb{R}^3} x_j B_{0i}(x) dx.
\ee
Finally we define
\be\no
\mathcal{M}_0=\{(u,B)\in [L^1(0,\infty; L^2(\mathbb{R}^3))]^6: \widetilde{\mathcal{A}}\ \text{is scalar and}\ \widetilde{\mathcal{C}}=\langle x, B_0\rangle\}.
\ee
Now we can state the lower bound results.
\bt\label{m3}(Lower bound).
{\it Let $(u_0,B_0)\in H\times H$ and $(u(x,t),B(x,t))$ be a weak solution of the Hall-MHD equations with initial datum $(u_0(x),B_0(x))$.
\begin{enumerate}[(1)]
  \item Assume that  the solution $(v,w)$ of (\ref{h}) satisfies
  \be\label{m301}
  c_1(t+1)^{-\alpha}\leq \|v(t)\|_2^2+ \|w(t)\|_2^2 \leq C_1(t+1)^{-\alpha}
  \ee
  for all $t\geq 0$, some constants $c, C>0$ and $0\leq \al<\f 52$. Then there exists $c_2, C_2>0$ such that
  \be\label{m302}
  c_2(t+1)^{-\alpha}\leq \|u(t)\|_2^2 +\|B(t)\|_2^2 \leq C_2(t+1)^{-\alpha}.
  \ee
  \item If $(u_0,B_0)\in [W_2\cap H\cap [L^1(\mathbb{R}^n)]^n]^2$ (so that $\widehat{u_0}(0)=0$ and $\widehat{B_0}(0)=0$), and $(u,B)\not\in \mathcal{M}_0$, then there exists $c_3, C_3>0$ such that
      \be\no
      c_3(t+1)^{-5/2}\leq \|u(\cdot,t)\|_2^2+ \|B(\cdot,t)\|_2^2 \leq C_3(t+1)^{-5/2}.
      \ee
\end{enumerate}
}
\et

In \cite{cdl}, the authors constructed the local and global in time strong solutions to the Hall-MHD (\ref{hmhd}). We record their results here as a Lemma.
\bl\label{m4}(Theorem 2.2 and 2.3 in \cite{cdl}).
{\it
Let $(u_0, B_0)\in H^{r}(\mathbb{R}^3)$ with $r>\f 52$, then there exists a unique strong solution $(u,B)\in L^{\infty}([0,T); H^r(\mathbb{R}^3))\cap Lip([0,T); H^{r-2}(\mathbb{R}^3))$ to the Hall-MHD equations (\ref{hmhd}) with $(u_0, B_0)$, where $T=T(\|u_0\|_{H^r}+\|B_0\|_{H^r})$. Moreover, there exists a positive constant $K_1(r)$, such that if $\|u_0\|_{H^r}+\|B_0\|_{H^r}\leq K_2(r)$, then $T=\infty$ and the solution becomes global.
}
\el
Here we are interested in the smoothing effect of these strong solutions. We will show that the local in time strong solution will become smooth after $t> 0$, indeed, it becomes analytic for strong solutions with small initial data. We expect that the strong solution with general initial data is also analytic at least local in time. Our method is based on Gevrey estimates developed in \cite{ft89},\cite{ot00} and the reference therein.

\bt\label{m5}
{\it
\begin{enumerate}[(1)]
  \item Let $(u,B)$ be a strong solution to the Hall-MHD equations (\ref{hmhd}), with initial value $(u_0,B_0)\in H^r(\mathbb{R}^3)$ with $r>\f 52$, then the local strong solution $(u,B)$ in Lemma \ref{m4} becomes smooth after $t>0$.
  \item There exists a constant $0<K_2\leq K_1$ such that if $\|(u_0,B_0)\|_{H^{r}}\leq K_2$, then the global strong solution $(u,B)$ in Lemma \ref{m4} becomes analytic after $t>0$ and the radius of analyticity grows like $\sqrt{t}$ in time. Furthermore, if there exists $\kappa_1>0$ and $\gamma\geq 0$ such that for any $t\geq 0$, there holds
      \be\label{m501}
      \|u(t)\|_2^2 +\|B(t)\|_2^2 \leq \f{\kappa_1}{(t+1)^{\gamma}}
      \ee
  Then there exists a constant $c_5=c_5(m,\kappa_1,\gamma)$ such that for every real number $m>0$
  \be\label{m502}
  \|\nabla^m u(t)\|_2^2+\|\nabla^m B(t)\|_2^2 \leq c_5 \f{1}{(t+1)^{\gamma+m}}.
  \ee
  \item If, in addition, $(u_0,B_0)\in L^1(\mathbb{R}^3)$ and there exists $\kappa_2,\kappa_3, \kappa_4$, which may depend on $(u_0,B_0)$, such that for $\forall \e>0$, there exists $t_1\geq0$, so that for all $t\geq t_1$,
  \be\label{m504}
  \|u(t)-v(t)\|_2^2+\|B(t)-w(t)\|_2^2 \leq \f{\e \kappa_2}{(t+1)^{\gamma}}
  \ee
  and for every $m\in \mathbb{N}$,
  \be\label{m505}
  \f{\kappa_3(m)}{(t+1)^{\ga+m}}\leq \|\nabla^m v(t)\|_2^2+\|\nabla^m w(t)\|_2^2\leq \f{\kappa_4(m)}{(t+1)^{\gamma}}.
  \ee
  Then there exists a positive constant $c_6=c_6(\kappa_2,\kappa_3,\kappa_4,\gamma, m)$ such that
  \be\label{m506}
  \|\nabla^m u(t)\|_2^2 +\|\nabla^m B(t)\|_2^2 \geq \f{c_6}{(1+t)^{\gamma+m}}.
  \ee
\end{enumerate}
}
\et

Here we remark that by Theorem \ref{m2}, we know that (\ref{m504}) will be satisfied when (\ref{m201}) holds for $0\leq \al<\f 52$. In the following, we will give the proof of Theorem \ref{m2}, \ref{m3}  and \ref{m5}. 

\section{Proof of Theorem \ref{m2} and \ref{m3} }

\subsection{Proof of Theorem \ref{m2}.} Let
\be\no
H(x,t) &=&(u\cdot\nabla) u- (B\cdot\nabla) B +\nabla \pi,\\\no
M(x,t) &=&(u\cdot\nabla) B- (B\cdot\nabla) u +\nabla\times ((\nabla\times B)\times B)\\\no
&=& (u\cdot\nabla) B- (B\cdot\nabla) u+ \nabla\times (\text{div }(B\otimes B)).
\ee
then Fourier transform of $(u,B)$ can be rewritten as
\be\label{3s8}
\begin{array}{lll}
\hat{u}(\xi,t) &=& e^{-t|\xi|^2} \widehat{u_0}(\xi) -\int_0^t e^{-(t-s)|\xi|^2} \hat{H}(\xi,s) ds,\\\no
\hat{B}(\xi,t) &=& e^{-t|\xi|^2} \widehat{B_0}(\xi) -\int_0^t e^{-(t-s)|\xi|^2} \hat{M}(\xi,s) ds.
\end{array}
\ee

Since $\nabla\cdot u=\nabla \cdot B=0$, applying the divergence operator to the first set of the Hall-MHD equations gives
\be\no
-\Delta \pi= \sum_{k,j=1}^n \f{\p^2}{\p x_k\p x_j} (u_k u_j-B_k B_j).
\ee
Hence
\be\label{3s9}
\hat{\pi}(\xi,t)= -\f{1}{|\xi|^2} \sum_{k,j} \xi_k\xi_j (\widehat{u_k u_j}- \widehat{B_k B_j}).
\ee
Then it follows that
\be\no
\hat{H}(\xi,t) &=& i\sum_{j} \xi_j(\widehat{u_j u}- \widehat{B_j B}) -i \sum_{k,j} \f{\xi_k\xi_j}{|\xi|^2}(\widehat{u_j u_k}- \widehat{B_j B_k}) \xi,\\\no
\hat{M}(\xi,t) &=& i\sum_{j} \xi_j (\widehat{u_j B}- \widehat{B_j u})- \xi\times(\xi_j\widehat{B_jB}).
\ee

Setting
\be\no
a_{kj}= \widehat{u_k u_j},\quad  b_{kj}= \widehat{B_k B_j},\quad c_{kj}= \widehat{u_j B_k}.
\ee
Introduce $A=[A_{kj}], C=[C_{kj}], \mu=[\mu_{kj}]$, where
\be\no
A_{kj}(\xi,t)= a_{kj}(\xi,t)- a_{jk}(\xi,t),\quad C_{kj}(\xi,t)= c_{kj}(\xi,t)- c_{jk}(\xi,t),\quad \mu_{kj}(\xi,t)=\f{\xi_k\xi_j}{|\xi|^2}.
\ee
Then
\be\no
\hat{H}(\xi,t) &=& i(I-\mu(\xi))A(\xi,t)\xi,\\\no
\hat{M}(\xi,t) &=& i C(\xi,t)\xi- \xi\times(\xi_j\widehat{B_jB}).
\ee
Since $I-\mu(\xi)$ is an orthogonal projection matrix for each $\xi \in \mathbb{R}^3\setminus \{0\}$, we get
\be\label{3s11}
|\hat{H}(\xi,t)| &\leq& C(\|u(t)\|_2^2+\|B(t)\|_2^2) |\xi|\\\label{3s13}
|\hat{M}(\xi,t)| &\leq& C(\|u(t)\|_2^2+\|B(t)\|_2^2) |\xi|+ C\|B(t)\|_2^2|\xi|^2.
\ee
By the energy estimate, we have
\be\no
\f{d}{dt}(\|u(t)\|_2^2+\|B(t)\|_2^2)= -(\|\nabla u(t)\|_2^2+\|\nabla B(t)\|_2^2).
\ee
Setting $E(t)=\|u(t)\|_2^2+\|B(t)\|_2^2$ and let $g(t)\geq 0$ for $t\geq 0$ (to be determined later) and $G(t)= \text{exp}\b(2\int_0^t g(s)^2 ds\b)$, so that $G'=2 g^2 G$. We have
\be\no
&\quad&\f{d}{dt}(G(t)E(t)) = 2G(t)(g(t)^2 E(t)-\|\nabla u(t)\|_2^2-\|\nabla B(t)\|_2^2)\\\label{3s16}
&\leq& 2g^2(t) G(t) \int_{|\xi|\leq g(t)} (|\hat{u}(\xi,t)|^2+ |\hat{B}(\xi,t)|^2) d\xi\\\no
&\leq& 2g^2(t) G(t)\int_{|\xi|\leq g(t)} |\hat{v}(\xi,t)|^2 + |\hat{w}(\xi,t)|^2+ \b(\int_0^t E(s) ds\b)^2(|\xi|^2+|\xi|^4) d\xi\\\no
&\leq& 2g^2(t)G(t)(\|v(t)\|_{2}^2+\|w(t)\|_2^2)+ g^2(t)(g^{5}(t)+ g^{7}(t))G(t)\b(\int_0^t E(s) ds\b)^2.
\ee
Suppose
\be\label{w31}
E(s) \leq C(1+s)^{-\beta},
\ee
with $\beta\geq 0$, we will use (\ref{3s16}) to improve the estimate on $E(t)$. By the energy inequality, (\ref{w31}) holds with $\beta=0$. Take $g^2(t)= \f{\ga}{2} (t+1)^{-1}$ with $\ga>\max\{1+\alpha, \f 32+ 2\beta\}$, and hence $G(t)=(t+1)^{\ga}$.
integrating (\ref{3s16}) over $[1,t]$, yielding
\be\no
E(t) (t+1)^{\gamma} \leq 2^{\ga} E(1) +c(t+1)^{\ga-\al}+ c(t+1)^{\ga-\f 12-2\beta}\quad \text{if $\beta<1$},
\ee
which improves the previous decay rate
\be\no
E(t) \leq C(t+1)^{-\bar{\beta}}
\ee
with $\bar{\beta}= \min\{\al, 2\beta+\f 12\}>\beta$. Start with this new exponent, and after finitely many iterations we conclude that
\be\no
E(t)\leq C(t+1)^{-\al}\quad \text{if $\al\leq 1$}.
\ee

If $\al>1$, after finitely many iterations we achieve $\bar{\beta}$ of the form $\bar{\beta}=1+\e$ with $\e>0$. Now
\be\no
\int_0^s E(r) dr\leq C,
\ee
which is independent of $s$, and by integrating (\ref{3s16}) for $\ga$ large, we obtain
\be\no
E(t) (t+1)^{\ga} \leq C + c(t+1)^{\ga-\al} + c(t+1)^{\ga-\f 52};
\ee
hence we finish the first part. For (2), since
\be\no
\p_t D_1=\Delta D_1 -H,\quad \p_t D_2 =\Delta D_2- M
\ee
and $D(0)=(0,0)$, we have
\be\no
\hat{D_1}(\xi,t)&=& -\int_0^t e^{-(t-s)|\xi|^2} \hat{H}(\xi,s) ds,\quad\hat{D_2}(\xi,t)=-\int_0^t e^{-(t-s)|\xi|^2} \hat{M}(\xi,s) ds.
\ee
Then by (\ref{3s11}) and (\ref{3s13}), we have
\be\label{3s19}
|\hat{D}(\xi,t)|&\leq& C(|\xi|+|\xi|^2) \int_0^t E(s) ds=: C(|\xi|+|\xi|^2) \Phi(t),\\\no
\f{d}{dt}\|D(t)\|_2^2 &=&- 2\|\nabla D(t)\|_2^2- 2\langle D_1, u\cdot\nabla v\rangle + 2\langle D_1, B\cdot\nabla w\rangle\\\label{3s14}
&\quad&- 2\langle D_2, u\cdot\nabla w\rangle+ 2\langle D_2, B\cdot\nabla v\rangle+ 2\langle \nabla\times D_2, (\nabla\times w)\times B\rangle\\\no
&\leq&- \|\nabla D(t)\|_2^2+ C\|\nabla w(t)\|_{\infty}^2\|B(t)\|_2^2+ C E(t)^{1/2} \|D(t)\|_2 (\|\nabla v(t)\|_{\infty}+\|\nabla w(t)\|_{\infty}).
\ee
Let $G(t)= \text{exp}\b(\int_0^t g^2(s) ds\b)$, then (\ref{3s14}) implies
\be\no
\f{d}{dt}(G(t)\|D(t)\|_2^2) &\leq& G(t)(g(t)^2 \|D(t)\|_2^2-\|\nabla D(t)\|_2^2)+ C G(t)\|\nabla\times w(t)\|_{L^{\infty}}^2\|B(t)\|_{2}^2\\\label{3s15}
&\quad&+ C G(t)E(t)^{1/2}\|D(t)\|(\|\nabla v\|_{\infty}+\|\nabla w\|_{\infty})\\\no
&\leq& g^2(t) G(t)\int_{|\xi|\leq g(t)} |\hat{D}(\xi,t)|^2 d\xi+ C G(t)\|\nabla\times w(t)\|_{L^{\infty}}^2\|B(t)\|_{2}^2\\\label{3s15}
&\quad&+ C G(t)E(t)^{1/2}\|D(t)\|(\|\nabla v\|_{\infty}+\|\nabla w\|_{\infty}).
\ee
By Lemma 2.6 in \cite{sss}, we have
\be\label{3s20}
\|\nabla v(t)\|_{\infty}^2 +\|\nabla w(t)\|_{\infty}^2\leq C(t+1)^{-\f 52-\alpha}.
\ee

Select $g(t)=\sqrt{\f{\gamma}{(t+1)}}$, then $G(t)=(1+t)^{\gamma}$, where $\gamma>\max\{\f 72,\f 52+2\alpha\}$. Together with (\ref{3s19}) and (\ref{3s20}), we obtain
\be\label{3s21}
\f{d}{dt}((t+1)^{\gamma}\|D(t)\|_2^2) &\leq& C \Phi(t)^2 (t+1)^{-\f 72+\gamma} + C \|D(t)\|_2 (t+1)^{-\f 54-\alpha+\gamma}\\\no
&\quad&\quad+ C(t+1)^{-\f 52-2\alpha+\gamma}.
\ee
Integrating over $[0,t]$, since $\Phi(t)$ is non-decreasing, we get
\be\no
\|D(t)\|_2^2 &\leq& C\Phi(t)^2 (t+1)^{-\f 52}+ C(t+1)^{-\gamma} \int_0^t (s+1)^{-\f 54-\alpha+\gamma} \|D(s)\|_2 ds+ C(t+1)^{-\f 32-2\alpha}.
\ee
Setting $Y(t)= \sup_{0\leq s\leq t} (s+1)^{5/4}\|D(s)\|_2$, the last inequality reduces
\be\no
Y(t)^2 \leq C\Phi(t)^2+ C(t+1)^{1-\alpha} Y(t)+ C (t+1)^{1-2\alpha},
\ee
thus
\be\no
Y(t)\leq C\Phi(t)+ C(t+1)^{1-\al}.
\ee
Since
\be\no
\Phi(t) =\int_0^t E(s) ds \leq \left\{\begin{array}{lll}
C(t+1)^{1-\al}\ \ \ &\text{if $0\leq \al<1$},\\
C\log(t+1)\ \ \ &\text{if $\alpha=1$},\\
C\ \ \ &\text{if $\alpha>1$},
\end{array}\right.
\ee
we have finished the proof of Theorem \ref{m2}.

\subsection{Proof of Theorem \ref{m3}.}The conclusions in (1) follows from Theorem \ref{m2} immediately. To prove (2), we follows the proof in \cite{sss} for the MHD case. Since $\widehat{u_0}(0)=\widehat{B_0}(0)=0$, it is well known that (\ref{m201}) holds with $\alpha=\f 52$, so $\|u(t)\|_2^2+ \|B(t)\|_2^2 \leq C(t+1)^{-\f 52}$.

We observe that $A,C$ are continuously differentiable in $\xi$ with bounded partial derivatives, although there is a new Hall term. This will be proved in Lemma \ref{3s161}. Then $A(\xi,t)=A(0,t)+ O_t(\xi)|\xi|, C(\xi,t)=C(0,t)+ O_t(\xi)|\xi|$. It follows that
\be\no
\hat{H}(\xi,t) &=& i(I-\mu(\xi))A(0,t)\xi+ O_t(\xi)|\xi|^2,\\\no
\hat{M}(\xi,t) &=& i C(0,t)\xi+ O_t(\xi)|\xi|^2.
\ee
We use this expansion in (\ref{3s8}) to get
\be\label{3s32}
\begin{array}{ll}
\hat{u}(\xi,t) &= P_1(\xi,t)\xi + O_t(\xi)|\xi|^2,\\
\hat{B}(\xi,t) &= P_2(\xi,t)\xi + O_t(\xi)|\xi|^2.
\end{array}
\ee
where
\be\no
P_1(\xi,t) &=& D_{\xi}\hat{u}_0(0)- i(I-\mu(\xi)) \mathcal{A}(t),\\\no
P_2(t) &=& D_{\xi}\hat{B}_0(0)- i\mathcal{C}(t).
\ee

Hence the expansion of $\hat{u}(\xi,t)$ and $\hat{B}(\xi,t)$ near $\xi=0$ is exactly same as the MHD case, then one can argue as in \cite{sss} to show that if $(u,B)\not\in \mathcal{M}_0$, then there exist $T_0>0, \rho>0$ such that either
\be\label{3s33}
\int_{\mathbb{S}^{n-1}} |P_1(\om,t)\om|^2 d\om\geq \rho
\ee
or
\be\label{3s35}
\int_{\mathbb{S}^{n-1}} |P_2(t)\om|^2 d\om\geq \rho
\ee
for $t\geq T_0$. Let $T\geq T_0$ (to be determined later) and let $(v(t),w(t))$ be the solution of the heat equation with initial datum $(v(0),B(0))=(u(T),B(T))$. In view of the representation (\ref{3s32}) (with $t=T$) for the initial datum of $(v,w)$, as was shown in Lemma 2.3 in \cite{sss} that there exists a constant $c>0$, such that
\be\no
\|v(t)\|_2^2 +\|w(t)\|_2^2 \geq c \rho t^{-5/2} +O(t^{-3}).
\ee
Now we will compare the solution $(u(t+T), B(t+T))$ with $(v(t), w(t))$. We set $D(t)=(D_1(t),D_2(t))=(u(t+T),B(t+T))-(v(t),w(t))$ so that $D$ satisfies
\be\no
\p_t D_1= \Delta D_1(t)- H(t+T),\quad \p_t D_2= \Delta D_2(t)- M(t+T)
\ee
and $D(0)=0$. As (\ref{3s21}), we have
\be\no
\f{d}{dt}(t^{\ga}\|D(t)\|_2^2) &\leq& t^{\ga-1}\int_{|\xi|\leq \ga/\sqrt{t}} |\hat{D}(\xi,t)|^2 d\xi+ C t^{\gamma}\|\nabla w(t)\|_{\infty}^2 \|B(t+T)\|_2^2\\\no
&\quad&+ C t^{\ga} E(t)^{1/2} \|D(t)\|_2 (\|\nabla v(t)\|_{\infty}+\|\nabla w(t)\|_{\infty}).
\ee
Note that
\be\no
|\hat{D}(\xi,t)|&\leq& C(|\xi|+|\xi|^2) \int_0^t (\|u(s+T)\|_2^2+\|B(s+T)\|_2^2) ds \leq C(|\xi|+|\xi|^2) T^{-3/2}.
\ee
This yields
\be\no
\f{d}{dt}(t^{\ga}\|D(t)\|_2^2) &\leq& C T^{-3} t^{\ga-7/2} + C_T t^{\ga-5}.
\ee
Taking $\ga>5$ and integrating over $[1,t]$, we finally obtain
\be\no
\|D(t)\|_2^2 \leq C T^{-3} t^{-5/2}+ C_T t^{-3}.
\ee
Taking $T$ large enough so that $C T^{-3}\leq \f 14 c\rho$, we get
\be\no
\|u(t+T)\|_2^2 + \|B(t+T)\|_2^2 &\geq& (\|v(t)\|_2+\|w(t)\|_2-\|D(t)\|_2)^2\\\no
&\geq&\f 12(\|v(t)\|_2^2+\|w(t)\|_2^2)-\|D(t)\|_2^2\\\no
&\geq&\f 14 c \rho t^{-5/2}+ O(t^{-3}).
\ee
It remains to establish the following lemma.
\bl\label{3s161}
{\it
Let$(u_0,B_0)$ belong to $[H^1(\mathbb{R}^3)\cap H\cap W_2]^2$. Suppose that $(u(t), B(t))$ are regular global solutions of the Hall-MHD equations with initial data $(u_0,B_0)$. Then for all $t\geq 0$,
\be\no
|\nabla_{\xi} a_{ij}(\xi,t)|+|\nabla_{\xi} b_{ij}(\xi,t)|+|\nabla_{\xi} c_{ij}(\xi,t)| &\leq& C(t),
\ee
where $a_{ij}=\widehat{u_i u_j}, b_{ij}=\widehat{B_i B_j}$ and $c_{ij}=\widehat{u_i B_j}$. Here $C(t)=(\|u_0\|_{W_2}^2+\|B_0\|_{W_2}^2)+c(1+\|u_0\|_2+\|B_0\|_2)^2(t+1)$.
}
\el

\bpf
Clearly,
\be\no
|\nabla_{\xi} a_{ij}|+|\nabla_{\xi} b_{ij}|+|\nabla_{\xi} c_{ij}| &\leq& C\int_{\mathbb{R}^3} |x| |u_i u_j| dx +
C\int_{\mathbb{R}^3} |x| |B_i B_j| dx + C\int_{\mathbb{R}^3} |x| |u| |B| dx \\\no
&\leq& C \int_{\mathbb{R}^3} |x|(|u|^2+|B|^2) dx.
\ee

It suffices to prove that
\be\no
\int_{\mathbb{R}^3} |x|(|u|^2+|B|^2) dx\leq C(t).
\ee

Dot-multiplying both sides of the first MHD equation with $|x|u$, of the second MHD equation with $|x|B$, adding and integrating over $\mathbb{R}^3$, we get after some integration by parts
\be\no
&\quad&\f{d}{dt} \int_{\mathbb{R}^3} |x|(|u|^2+|B|^2) dx =- \int_{\mathbb{R}^3}|x|(|\nabla u|^2+|\nabla B|^2) dx+ \int_{\mathbb{R}^3} \f{|u|^2+|B|^2}{|x|} dx\\\no
&+&\f 12\int_{\mathbb{R}^3} \f{(x\cdot u)|u|^2}{|x|} dx- \int_{\mathbb{R}^3} \f{(x\cdot B)(u\cdot B)}{|x|} dx+\f 12\int_{\mathbb{R}^3}\f{(x\cdot u)|B|^2}{|x|} dx+\int_{\mathbb{R}^3}\f{1}{|x|}(x\cdot u)p dx\\\no
&-& \int_{\mathbb{R}^3}\b(\f{x}{|x|}\times B\b)\cdot ((\nabla \times B)\times B) dx\\\no
&\leq&\int_{\mathbb{R}^3} \f{1}{|x|}(|u|^2+|B|^2) dx+\f 12 \int_{\mathbb{R}^3} |u|^3 dx +2\int_{\mathbb{R}^3} |u||B|^2 dx+ \int_{\mathbb{R}^3} |u||p| dx+\int_{\mathbb{R}^3} |\nabla B| |B|^2 dx\\\no
&=:& I_1 +I_2 +I_3 + I_4+ I_5.
\ee
We estimate $I_i, i=1,\cdots, 5$ as follows.
\be\no
I_2+I_3+I_4 &\leq& C\|u\|_2 (\|u\|_4^2 +\|B\|_4^2)\\\no
|I_5|&\leq& C\|\nabla B\|_2\|B\|_{4}^2\leq C\|\nabla B\|_2 \|B\|_2\|B\|_6\leq C\|\nabla B\|_2^2 \|B\|_2,\\\no
|I_1|&=& \int_{|x|\leq 1} \f{1}{|x|} (|u|^2+|B|^2) dx+ \int_{|x|>1} \f{1}{|x|} (|u|^2+|B|^2) dx\\\no
&\leq& C(\|u\|_6^2+\|B\|_6^2)+ (\|u\|_2^2+\|B\|_2^2).
\ee

By the Sobolev embedding theorem, we have
\be\no
\f{d}{dt}\int_{\mathbb{R}^3}|x| (|u|^2+|B|^2) dx\leq C (\|u\|_{H^1}^2+\|B\|_{H^1}^2)(1+\|u\|_2+\|B\|_2).
\ee
This yields
\be\no
\int_{\mathbb{R}^3}|x| (|u|^2+|B|^2) dx &\leq&(\|u_0\|_{W_2}^2+\|B_0\|_{W_2}^2)\\\no
&\quad&+ C(1+\|u_0\|_2+\|B_0\|_2) \int_0^t (\|u(s)\|_{H^1}^2+\|B(s)\|_{H^1}^2) ds\\\no
&\leq& (\|u_0\|_{W_2}^2+\|B_0\|_{W_2}^2)+C(1+\|u_0\|_2+\|B_0\|_2)^2(t+1).
\ee

\epf

\section{Proof of Theorem \ref{m5}}

\subsection{Smoothing effects of the local strong solution}

We start with the proof of (1) in Theorem \ref{m5}. We need the following energy estimate, which is slightly different from those in \cite{cdl}: for any integer $m>\f 52$, there exists a constant $C(m)>0$, such that
\be\label{m5011}
&\quad&\f 12 \f{d}{dt}(\|u\|_{H^m}^2+\|B\|_{H^m}^2) + \|\nabla u\|_{H^m}^2+ \|\nabla B\|_{H^m}^2\\\no
&\leq& C(m)(\|u\|_{H^m}^2+\|B\|_{H^m}^2)(\|\nabla u\|_{H^m}+\|\nabla B\|_{H^m}).
\ee
The inequality (\ref{m5011}) follows from the simple energy estimates. For any multi-index $\alpha\in \mathbb{N}_0^3$ with $m=|\al|$, from the Hall-MHD equations (\ref{hmhd}), we have
\be\no
&\quad&\f{1}{2} \f{d}{dt}(\|u\|_{H^m}^2+\|B\|_{H^m}^2)+ \|\nabla u\|_{H^m}^2+\|\nabla B\|_{H^m}^2\\\no
&=&-\sum_{0\leq |\alpha|\leq m} \int \nabla^{\alpha} (u\cdot\nabla u)\cdot \nabla^{\alpha} u dx+ \sum_{0\leq |\alpha|\leq m} \int \nabla^{\alpha} (B\cdot\nabla B)\cdot \nabla^{\alpha} u + \nabla^{\alpha}(B\cdot u)\cdot \nabla^{\alpha} B dx\\\no
&\quad&-\sum_{0\leq |\alpha|\leq m} \int \nabla^{\alpha} (u\cdot\nabla B)\cdot \nabla^{\alpha} B dx-\sum_{0\leq |\alpha|\leq m} \int \nabla^{\alpha}((\nabla\times B)\times B)\cdot \nabla^{\alpha}(\nabla\times B) dx\\\no
&=:& I +II + III +IV.
\ee
By using the calculus inequality:
\be\no
\sum_{|\alpha|\leq m} \|\nabla^{\al}(fg)- (\nabla^{\alpha} f) g\|_2 &\leq& C(m)(\|f\|_{H^{m-1}}\|\nabla g\|_{\infty}+\|f\|_{\infty}\|g\|_{H^m}),
\ee
we can bounded $I$ to $IV$ as follows, which immediately yields (\ref{m5011}) by Sobolev embedding theorem.
\be\no
|I|&=& \b|-\sum_{0\leq |\alpha|\leq m} \int [\nabla^{\alpha}(u\cdot\nabla u)-u\cdot\nabla \nabla^{\alpha}u] \cdot D^{\alpha} u dx\b|\leq C\|u\|_{H^m}^2 \|\nabla u\|_{\infty},\\\no
|II|&=&\b|\sum_{0\leq |\alpha|\leq m} \int [\nabla^{\alpha}(B\cdot \nabla u)-B\cdot\nabla \nabla^{\alpha} u]\cdot D^{\alpha} B\\\no
&\quad&\quad+[\nabla^{\alpha}(B\cdot\nabla B)-B\cdot\nabla \nabla^{\alpha} B]\cdot \nabla^{\alpha} u dx\b|\\\no
&\leq& C\|u\|_{H^m}\|B\|_{H^m}\|\nabla B\|_{\infty}+ C\|\nabla u\|_{\infty}\|B\|_{H^m}^2,\\\no
|III|&=&\b|-\sum_{0\leq |\alpha|\leq m}\int [\nabla^{\alpha}(u\cdot\nabla B)-u\cdot \nabla \nabla^{\al} B]\cdot \nabla^{\alpha} B dx\b|\\\no
&\leq& C(\|\nabla B\|_{H^{m-1}}\|\nabla u\|_{\infty}+\|\nabla B\|_{\infty}\|u\|_{H^m})\|B\|_{H^m}.
\ee
\be\no
|IV|&=&\b|-\sum_{0\leq |\alpha|\leq m}\int [\nabla^{\alpha}((\nabla\times B)\times B)-\nabla^{\alpha}(\nabla\times B)\times B]\cdot \nabla^{\alpha}(\nabla\times B) dx\b|\\\no
&\leq& C(\|\nabla\times B\|_{H^{m-1}}\|\nabla B\|_{\infty}+\|\nabla \times B\|_{\infty}\|B\|_{H^m})\|\nabla B\|_{H^m}\\\no
&\leq& C\|B\|_{H^m}\|\nabla B\|_{\infty}\|\nabla B\|_{H^m}.
\ee
For $\forall t_1\in (0,T)$, we want to show that $\|u(t_1)\|_{H^m}+\|B(t_1)\|_{H^m}<\infty$ for any positive integer $m\in \mathbb{N}$. By the local existence theorem in \cite{cdl}, there exists $C_1,C_2>0$ which may depend on $t_1$ such that
\be\label{m5012}
\begin{array}{ll}
\|u(t)\|_{H^r}+\|B(t)\|_{H^r}\leq C_1(\|u_0\|_{H^r}+\|B_0\|_{H^r}),\\
\int_0^{t_1} \|\nabla u(s)\|_{H^r}+\|\nabla B(s)\|_{H^r} ds \leq C_2(\|u_0\|_{H^r}+\|B_0\|_{H^r}).
\end{array}
\ee
Without loss of generality, we assume that $r\in \mathbb{N}$. By (\ref{m5012}), there exists $t_2\in (0, t_1)$ such that $\|\nabla u(t_2)\|_{H^r}+\|\nabla B(t_2)\|_{H^r}<\infty$, hence $\|u(t_2)\|_{H^{r+1}}+\|B(t_2)\|_{H^{r+1}}<\infty$. We integrate (\ref{m5011}) with $m=r+1$ over $[t_2,t_1]$, yielding that for $\forall t\in [t_2,t_1]$
\be\label{m5013}
\|u(t)\|_{H^{r+1}}^2+\|B(t)\|_{H^{r+1}}^2 \leq e^{\int_{t_2}^{t} (\|\nabla u(s)\|_{H^r}+\|\nabla B(s)\|_{H^r}) ds} (\|u(t_2)\|_{H^{r+1}}^2+\|B(t_2)\|_{H^{r+1}})^2<\infty.
\ee
Especially, we have $\|u(t_1)\|_{H^{r+1}}+\|B(t_1)\|_{H^{r+1}}<\infty$. With (\ref{m5013}), (\ref{m5011}) also produces
\be\label{m5014}
\int_{t_2}^{t_1} \|\nabla u(s)\|_{H^{r+1}}+\|\nabla B(s)\|_{H^{r+1}} ds <\infty.
\ee
By (\ref{m5013}) and (\ref{m5014}), we can find another $t_3\in (t_2,t_1)$ such that $\|u(t_3)\|_{H^{r+2}}+\|B(t_3)\|_{H^{r+2}}<\infty$ and then argue as previous to show that $\|u(t_1)\|_{H^{r+2}}+\|B(t_1)\|_{H^{r+2}}<\infty$. Continuing this process, we finish the proof.

\subsection{Analyticity of small strong solutions and upper bound}

We will use Gevrey estimates to show the analyticity of the strong solutions to the Hall-MHD. Setting $\Lambda=(-\Delta)^{1/2}$, for $\tau\geq 0$, we introduce the spaces
\be\no
\mathcal{D}(e^{\tau \Lambda}; H^r)= \{w\in H^r(\mathbb{R}^3): e^{\tau\Lambda} w\in H^r(\mathbb{R}^3)\}.
\ee
As shown in \cite{ot00}, for every $w\in \mathcal{D}(e^{\tau \Lambda}; H^r)$ with $\tau>0, r>0$, then for $\forall x\in \mathbb{R}^3$ and every multi-index $\alpha\in \mathbb{N}_0^3$, there exists $M$ and $\rho=\tau/\sqrt{3}$, such that
\be\no
|\p^{\alpha} w(x)|\leq M\f{\alpha !}{\rho^{|\alpha|}}.
\ee
That is, $w$ is analytic with radius $\tau/\sqrt{3}$ on the whole of $\mathbb{R}^3$. In the following, we only need to show that the strong solution $(u,B)$ belongs to $\mathcal{D}(e^{\tau \Lambda}; H^r)$ with $\tau>0$. First, we need the following Lemmas, which was proved in \cite{ot00}.

\bl\label{m11}
{\it
Let $\tau\geq 0, r>3/2$, and $s<3/2$. Then there exists a constant $C=C(r,s)$ such that any two functions $v$ and $w$ in $\mathcal{D}(e^{\tau \Lambda}; H^r)$ satisfy the inequality
\be\label{m111}
\|\Lambda^r e^{\tau \Lambda}(vw)\|_{2} &\leq& C(r,s) (\|\Lambda^s e^{\tau \Lambda}w\|_{H^{r-s}}\|\Lambda^r e^{\tau \Lambda} v\|_{2}+\|\Lambda^s e^{\tau \Lambda}v\|_{H^{r-s}}\|\Lambda^r e^{\tau \Lambda} w\|_{2}).
\ee
}
\el

\bl\label{m12}
{\it
The following inequalities hold:
\be\label{m121}
\|\Lambda^r e^{\tau \Lambda} u\|_2^2 &\leq& 2\|\Lambda^r u\|_2^2 + 2 \tau^2 \|\Lambda^{r+1} e^{\tau \Lambda} u\|_2^2,\quad \text{for all $r\geq 0$ and $\tau\geq 0$};\\\label{m122}
\|\Lambda^p e^{\tau \Lambda} u\|_2^2 &\leq& e\|\Lambda^p u\|_2^2 + (2\tau)^{2q} \|\Lambda^{p+q} e^{\tau \Lambda} u\|_2^2,\quad \text{for all nonnegative $p,q$ and $\tau$;}\\\label{m123}
\|\Lambda^q u\|_2^2 &\leq& c(p,q) \tau^{p-2q} \|u\|_2 \|\Lambda^p e^{\tau \Lambda} u\|_2, \quad \text{for $2q\geq p\geq 0$ and $\tau> 0$}.
\ee
}
\el

Now we start to prove the conclusion (2) in Theorem \ref{m5}. Setting
\be\no
\begin{array}{ll}
J_r&= \|\Lambda^r u\|_{2}^2,\quad H_r=\|\Lambda^r B\|_{2}^2,\quad G_r= \|\Lambda^r e^{\tau \Lambda} u\|_{2}^2,\quad K_r= \|\Lambda^r e^{\tau \Lambda} B\|_{2}^2,\\
N_r&= H_r+ J_r,\quad M_r=G_r+K_r.
\end{array}
\ee
First, we show the local in time analyticity by choosing $\tau=\tau(t)$. In the following, we assume $(u,B)$ are the global strong solution to the Hall-MHD with small initial data $(u_0,B_0)$ in Lemma \ref{m4}. Then $N_r(t)\leq C(r) N_r(0)$. And
\be\no
&\quad&\f 12 \f{d}{dt} M_r=\tau'(t)\int_{\mathbb{R}^3} \Lambda^{r+1} e^{\tau \Lambda} u \cdot \Lambda^r e^{\tau \Lambda} u+ \Lambda^{r+1} e^{\tau \Lambda} B\cdot \Lambda^r e^{\tau \Lambda} B\\\no
&\quad&\quad\quad +\int_{\mathbb{R}^3} \Lambda^r e^{\tau \Lambda} u\cdot \Lambda^r e^{\tau \Lambda} u_t+ \Lambda^r e^{\tau \Lambda} B\cdot \Lambda^r e^{\tau \Lambda} B_t\\\no
&=&\tau'(t)M_{r+1/2}-M_{r+1}- \int \Lambda^r e^{\tau \Lambda}(u\cdot\nabla u)\cdot \Lambda^r e^{\tau \Lambda} u+\int \Lambda^r e^{\tau \Lambda}(B\cdot\nabla B) \Lambda^r e^{\tau \Lambda} u\\\no
&\quad&-\int \Lambda^r e^{\tau \Lambda}(u\cdot\nabla B) \Lambda^r e^{\tau \Lambda} B+\int \Lambda^r e^{\tau \Lambda}(B\cdot\nabla u)\cdot \Lambda^r e^{\tau \Lambda} B\\\no
&\quad&-\int\nabla\times (\Lambda^r e^{\tau \Lambda})((\nabla\times B)\times B)\cdot \Lambda^r e^{\tau \Lambda} B\\\label{m18}
&\eqqcolon &\tau'(t)M_{r+1/2}-M_{r+1} +\sum_{i=1}^5 I_i.
\ee
By Lemma \ref{m11}, Lemma \ref{m12}, with $r>\f 32$ and $s<\f 32$, we have
\be\no
|I_1|&\leq& \|\Lambda^r e^{\tau \Lambda}(u\cdot\nabla u)\|_{2}\|\Lambda^r e^{\tau \Lambda} u\|_{2} \\\no
&\leq& c(\|\Lambda^s e^{\tau \Lambda}\nabla u\|_{H^{r-s}}\|\Lambda^r e^{\tau \Lambda} u\|_{2}+ \|\Lambda^s e^{\tau \Lambda} u\|_{H^{r-s}}\|\Lambda^r e^{\tau \Lambda} \nabla u\|_{2}) G_r^{1/2}\\\no
&\leq& c(G_s^{1/2}G_r^{1/2} G_{r+1}^{1/2}+G_{s+1}^{1/2} G_r+ G_r G_{r+1}^{1/2})\\\no
&\leq& c[(J_s^{1/2}+\tau^{r-s}G_r^{1/2})G_r^{1/2}G_{r+1}^{1/2}+(J_{s+1}^{1/2}+\tau^{r-s} G_{r+1}^{1/2})G_r+G_r G_{r+1}^{1/2}]\\\no
&\leq& c (J_s M_r+ J_{s+1}^{1/2}) M_r+ c (1+\tau^{r-s})^2 M_r^2+ \f{1}{100}M_{r+1}.
\ee
Similarly, we have
\be\no
\sum_{i=2}^4|I_i|&\leq&c (1+H_s+J_s+H_{s+1}^{1/2}+ J_{s+1}^{1/2}+H_{s+1})M_r+c(1+\tau^{r-s})^2 M_{r}^2+ \f{1}{10}M_{r+1}\\\no
&\leq& c(N_r) M_r+ c(1+\tau^{r-s})^2 M_r^2 +\f{1}{10} M_r,\\\no
|I_5| &=& \b|\int (\Lambda^r e^{\tau \Lambda})((\nabla\times b)\times b)\cdot (\nabla\times \Lambda^r e^{\tau \Lambda} b) dx\b|\\\no
&\leq& c\b((K_{s+1}^{1/2}+K_{r+1}^{1/2})K_r^{1/2}+(K_s^{1/2}+K_r^{1/2})K_{r+1}^{1/2}\b) K_{r+1}^{1/2}\\\no
&\leq& c H_{s+1} M_r +c (H_{s}^{1/2}+(1+\tau^{r-s})M_r^{1/2}) M_{r+1}+\f{1}{10} M_{r+1}.
\ee
where $c(N_r)$ is a smooth function of $N_r$, the function $c(N_r)$ may change in different line. Since we are considering the local in time analyticity, we restrict ourself on $t\in [0,1]$, and select $\tau=\tau(t)=t$, then $M_r(0)= N_r(0)$ and
\be\no
\f12 \f{d}{dt} M_r &\leq& c M_{r}^{1/2} M_{r+1}^{1/2}- \f 34 M_{r+1}+ c(N_r) M_r+ c M_r^2+ c(N_r^{1/2}+ M_r^{1/2}) M_{r+1}\\\label{m20}
&\leq& c(N_r) M_r + c M_r^2-\b(\f 12- c(N_r^{1/2}+ M_r^{1/2})\b) M_{r+1}.
\ee
Suppose $c(N_r^{1/2}(0)+ M_r^{1/2}(0))= 2c(\|u_0\|_{H^r}+\|B_0\|_{H^r})<\f 14$, then in short time, (\ref{m20}) reduces to
\be\no
\f{d}{dt} M_r \leq c(N_r(0)) M_r + c M_r^2.
\ee
Hence we can choose $N_r(0)$ small enough so that in short time $[0,\sigma]$, such that for $\forall t\in [0,\sigma]$
\be\label{m200}
M_r(t)\leq 2 M_r(0)= 2(\|u_0\|_{H^r}+\|B_0\|_{H^r})^2,
\ee
which can also guarantee that $c(N_r^{1/2}(t)+ M_r^{1/2}(t))<\f 12$. So we have showed that $M_r(t)$ is finite in some time interval $[0,\sigma]$.

Now we will refine our estimate to show that $M_r$ is finite at any time. Without loss of generality, we assume that $(u_0, B_0)\in \mathcal{D}(e^{\eta\Lambda}; H^r)$ for some $\eta>0$. The point is to explore the dissipation term $M_{r+1}$. Indeed, by Lemma \ref{m12}, we have $\f{M_r-2 N_r}{2\tau^2}\leq M_{r+1}$ and $N_r\leq \f{1}{\tau^{2r}} N_0+ \f{1}{8} M_r$. Hence (\ref{m18}) reduces
\be\no
\f 12 \f{d}{dt} M_r &\leq& \f{1}{2}\f{\tau'(t)}{\tau} M_r+\f12 \tau \tau'M_{r+1}- M_{r+1} + \sum_{i=1}^5 I_i\\\no
&\leq&\f{1}{2}\f{\tau'(t)}{\tau} M_r+\f12 \tau \tau'M_{r+1}-\f{M_r-2 N_r}{4\tau^2}-\f 12 M_{r+1}+ \sum_{i=1}^5 I_i\\\label{m21}
&\leq&\b(\f 12\f{\tau'}{\tau}-\f{1}{8\tau^2}\b) M_r+\b(\f 12 \tau'\tau-\f 18\b)M_{r+1}-\f{1}{8\tau^2} M_r\\\no
&\quad&\quad+\f{1}{2\tau^2}(\f{1}{\tau^{2r}}N_0+\f{1}{8}M_r)-\f{3}{8}M_{r+1}+ \sum_{i=1}^5 I_i.
\ee
We choose $\tau(t)=\sqrt{\tau_0^2+\alpha t}$, where $\tau_0>0$ and $0<\alpha\leq \f 12$ will be determined later. The point is $\f 12\tau\tau'=\f{\alpha}{4}\leq \f 18$, then (\ref{m21}) will reduce to
\be\label{m22}
\f{1}{2}\f{d}{dt}M_r &\leq&-\f{1}{16\tau^2} M_r-\f{3}{8}M_{r+1}+\f{1}{2\tau^{2r+2}}N_0+ \sum_{i=1}^5 I_i.
\ee
We will refine our estimates on $I_i, i=1,\cdots,5$ by using the ``good" term $M_{r+1}$. We will also replace $J_s,H_s, J_{s+1}, H_{s+1}$ by $J_0,H_0$ and $M_r$, since by assumptions, we have got the decay rates for $H_0, J_0$. By Lemma \ref{m12}, we have
\be\no
J_{s} &\leq& c(r,s) \tau^{r-2s} J_0^{1/2} G_r^{1/2},\ \ J_{s+1}\leq c(r,s)\tau^{r-2s-2} J_0^{1/2} G_r^{1/2},\\\no
H_{s} &\leq& c(r,s) \tau^{r-2s} H_0^{1/2} K_r^{1/2},\ \ H_{s+1}\leq c(r,s)\tau^{r-2s-2} H_0^{1/2} K_r^{1/2}.
\ee
Then the bounds for $I_i, i=1,\cdots, 5$ are followed in order.
\be\no
|I_1|&\leq& c_1 J_s^{1/2} G_r^{1/2} G_{r+1}^{1/2}+ c_1J_{s+1}^{1/2} G_r+ c_1(1+\tau^{r-s}) G_r G_{r+1}^{1/2}\\\no
&\leq& c_1 \tau^{r/2-s} J_0^{1/4} G_r^{3/4} G_{r+1}^{1/2}+ c_1 \tau^{r/2-s-1} J_0^{1/4} G_r^{5/4}+ c_1(1+\tau^{r-s}) G_r G_{r+1}^{1/2}\\\no
&\leq& c\tau^{r-2s} J_0^{1/2} M_r^{3/2}+ c \tau^{r/2-s-1} J_0^{1/4} M_r^{5/4}+ c(1+\tau^{2r-2s}) M_r^2+\f{3}{80} G_{r+1},
\ee
Similarly, we also have
\be\no
\sum_{i=2}^5 |I_i|&\leq& c\tau^{r-2s} N_0^{1/2} M_r^{3/2}+ c \tau^{r/2-s-1} N_0^{1/4} M_r^{5/4}+ c(1+\tau^{2r-2s}) M_r^2+\f{3}{80} K_{r+1}.
\ee

Back to (\ref{m22}), we obtain
\be\no
\f{1}{2}\f{d}{dt} M_r &\leq& -\f{1}{16\tau^2}M_r-\f{3}{16} M_{r+1} + \f{c}{\tau^{2r+2}}N_0\\\no
&\quad&+c\b[\tau^{r-2s}N_0^{1/2} M_r^{1/2}+ \tau^{r/2-s-1}N_0^{1/4} M_r^{1/4}+ (1+\tau^{2r-2s}) M_r\b] M_r\\\label{m23}
&\quad&+c[\tau^{r/2-s} H_0^{1/4} M_r^{1/2}+(1+\tau^{r-s}) M_r^{1/2}] M_{r+1}.
\ee
For our purpose, we want to choose initial data $(u_0,B_0)$ small enough, such that
\be\label{m31}
&&g_1(\tau):=c[\tau^{r/2-s} H_0^{1/4} M_r^{1/2}+(1+\tau^{r-s}) M_r^{1/2}]<\f{3}{16},\\\label{m32}
&&g_2(\tau):=c\b[\tau^{r-2s}N_0^{1/2}M_r^{1/2}+ \tau^{r/2-s-1} N_0^{1/4} M_r^{1/4}+ (1+\tau^{2r-2s})M_r\b]<\f{1}{32\tau^2}.
\ee
By (\ref{m200}), we have (\ref{m31}) holds in $[0,\sigma]$. Note that at $t=0$, $\|\Lambda e^{\eta_0\Lambda}(u_0,B_0)\|_{2}^2$ for $0\leq \eta_0\leq \eta$ is bounded by $\|\Lambda^r(u_0, B_0)\|_{2}^2<\infty$ and $\|\Lambda e^{\eta\Lambda}(u_0,B_0)\|_{2}^2<\infty$. By choosing $s\in [r/2, r/2+1)$, the powers of $\tau$ in $g_2$ is less than $2$, so that $\f{1}{32\tau^2}$ diverges faster as $\tau\to 0$. Then we can choose $\tau(0)=\tau_0\in (0,\eta]$ small enough that (\ref{m32}) is satisfied at $t=0$. Moreover, the differential inequality (\ref{m23}) admits a local smooth solution, then (\ref{m32}) is satisfied near $t=0$. The restriction $s\in [r/2, r/2+1)$ and $s<\f 32, r>\f 32$ will require $r<3$. Here for convenience, we choose $s=r/2=\f{11}{8}$. However, this is not a serious restriction, since the initial data $(u_0,B_0)$ in $H^r(\mathbb{R}^3)$ with $r>3$ is automatically in $H^{11/4}(\mathbb{R}^3)$.

So as long as (\ref{m31})-(\ref{m32}) are satisfied, (\ref{m23}) is reduced to
\be\no
\f{d}{dt} M_r &\leq& -\f{1}{16\tau^2} M_r+\f{c}{\tau^{2(r+1)}} N_0.
\ee
Note that if we choose $\alpha\leq \tau_0^2$, then $(1+t)^{-1/2}\leq (1+\f{\al}{\tau_0^2}t)^{-1/2}=\tau_0/\tau$, so that $H_0+J_0\leq \kappa_1(\tau_0/\tau)^{2\gamma}$. Then we conclude that
\be\no
\f{d}{dt}M_r &\leq&-\f{1}{16\tau^2} M_r+\f{c \kappa_1}{\tau^{2(\gamma+r+1)}}.
\ee
Multiplying the corresponding integral factor, it produces
\be\no
\f{d}{dt}(\tau^{\f{1}{8\alpha}}M_r) \leq c \tau^{2(\f{1}{16\alpha}-\gamma-r-1)}.
\ee
Choosing $\alpha$ small enough such that $1>16\al(\gamma+r)$, we obtain
\be\no
M_r(t)&\leq& \b(M_{r}(0)-\f{16 c}{1-16\al(\gamma+r)}\f{1}{\tau_0^{2(\gamma+r)}}\b)\b(\f{\tau_0^2}{\tau^2}\b)^{\f{1}{16\alpha}}\\\no
&\quad&\quad\quad +\f{16 c\kappa_1}{1-16\al(\ga+r)}\f{1}{\tau^{2(\gamma+r)}}\\\label{m35}
&\leq&c \kappa_1 \tau^{-2(\gamma+r)},
\ee
if we choose
\be\no
M_r(0)\leq \f{16 c}{1-16\al(\gamma+r)}\f{1}{\tau_0^{2(\gamma+r)}}.
\ee
Estimate (\ref{m35}) and the choice of $s=r/2=\f{11}{8}$ shows that
\be\no
g_1(\tau)&\leq& c\kappa_1(\tau^{-\ga-s-r/2}+\tau^{-\ga-s})\leq \kappa_1 \tau_0^{-s}<\f{3}{16},\\\no
g_2(\tau)&\leq& c\kappa_1(\tau^{-\ga-2s}+ \tau^{-s-1-\ga/2}+\tau^{-2\ga-2s})\leq c\kappa_1 \tau^{-s-1}<\f{1}{32\tau^2}
\ee
if $\kappa_1$ is small enough, which is guaranteed by a small initial data.

Finally, we can conclude that there exists a $K_2>0$ such that $\|(u_0,B_0)\|_{H^r}\leq K_2$, then for all $t\geq 0$
\be\no
M_r(t)\leq c\kappa_1 \tau^{-2(\gamma+r)}.
\ee
By Lemma \ref{m12}, we have
\be\no
\|\La^m u(t)\|_2^2+\|\La^m B(t)\|_2^2 &\leq & c(m,r) \tau^{r-2m} N_0^{1/2} M_r^{1/2} \\\no
&\leq& c(m,r) \tau^{-2(\gamma+m)}.
\ee

\subsection{The lower bound for higher order derivatives} In the proof of (3) in Theorem \ref{m5}, we will do the Gevrey estimates for the difference $D=(u-v,B-W)$ between the Hall-MHD and heat system. Setting
\be\no
\mathcal{M}_r(t)=\|\Lambda^r e^{\tau \La} D_1(t)\|_2^2+\|\La^r e^{\tau\La} D_2(t)\|_2^2.
\ee
Same as previous, we derive
\be\no
\f{d}{dt} \mathcal{M}_r &=& \tau'(t) \mathcal{M}_{r+1/2}- \mathcal{M}_{r+1}- \int \Lambda^r e^{\tau \Lambda}(u\cdot\nabla u)\cdot \Lambda^r e^{\tau \Lambda} D_1+\int \Lambda^r e^{\tau \Lambda}(B\cdot\nabla B) \Lambda^r e^{\tau \Lambda} D_1\\\no
&\quad&-\int \Lambda^r e^{\tau \Lambda}(u\cdot\nabla B) \Lambda^r e^{\tau \Lambda} D_2+\int \Lambda^r e^{\tau \Lambda}(B\cdot\nabla u)\cdot \Lambda^r e^{\tau \Lambda} D_2\\\no
&\quad&-\int\nabla\times (\Lambda^r e^{\tau \Lambda})((\nabla\times B)\times B)\cdot \Lambda^r e^{\tau \Lambda} D_2\\\no
&\eqqcolon &\tau'(t) \mathcal{M}_{r+1/2}- \mathcal{M}_{r+1}+ \sum_{i=1}^5 I_i',
\ee
where
\be\no
|I_1'|&\leq&C(r,s)(G_s^{1/2}G_{r+1}^{1/2}+ G_{s+1}^{1/2} G_r^{1/2}+ G_r^{1/2} G_{r+1}^{1/2}) \mathcal{M}_r^{1/2}\\\no
&\leq& C(r,s)(M_s M_{r+1}+ M_{s+1} M_r+ M_r M_{r+1})+ \f{1}{200} \mathcal{M}_r,\\\no
\sum_{i=2}^4|I_i'|&\leq& C(r,s)(M_s M_{r+1}+ M_{s+1} M_r+ M_r M_{r+1})+ \f{1}{200} \mathcal{M}_r.
\ee
\be\no
|I_5'|&=&\b|\int \La^r e^{\tau \La}((\nabla\times B)\times B)\cdot (\nabla \times \La^r e^{\tau \La} D_2) dx\b|\\\no
&\leq& C(r,s)(K_s^{1/2}K_{r+1}^{1/2}+ K_{s+1}^{1/2} K_r^{1/2}+ K_r^{1/2} K_{r+1}^{1/2}) \mathcal{M}_{r+1}^{1/2}\\\no
&\leq& C(r,s)(M_s M_{r+1}+ M_{s+1} M_r+ M_r M_{r+1})+ \f{1}{200} \mathcal{M}_{r+1}.
\ee
Then
\be\no
\f 12 \f{d}{dt}\mathcal{M}_r &\leq& \tau'(t)\mathcal{M}_{r+1/2}- \mathcal{M}_{r+1}+ C(M_s M_{r+1}+M_{s+1}M_r+ M_r M_{r+1})\\\no
&\quad&\quad+\f{1}{40} \mathcal{M}_r+ \f{1}{200}\mathcal{M}_{r+1}\\\no
&\leq& -\f{1}{16\tau^2} M_r+\f{1}{2\tau^2}(H_0+J_0) +C(M_s M_{r+1}+M_{s+1}M_r+ M_r M_{r+1})\\\no
&\leq& -\f{1}{16\tau^2} M_r+\f{c_8\e}{\tau^{2(\gamma+r+1)}}+ O(\tau^{-4\ga-2s-2r-2}).
\ee
By integrating as above, we finally get
\be\no
\mathcal{M}_r(t) \leq \f{c_9\e}{\tau^{2(\gamma+r)}}+ O(\tau^{-4\gamma-2s-2r})+ O(\tau^{-32\alpha}),
\ee
which implies that
\be\no
\|\La^m D(t)\|_2^2 \leq \f{\e c(m,r)}{\tau^{2(\gamma+m)}}.
\ee
For a given $m$, we choose $\e$ small enough so that $\kappa_3(m)> c(m,r)\e$, whence the triangle inequality implies the required lower bound.

{\bf Acknowledgements.} The author would like to thank Prof. Dongho Chae for his interest and stimulating discussions and the kind hospitality during the author's visit at Chung Ang University on December, 2014. Special thanks also go to Prof. Jiahong Wu for his interest and discussion.

\bibliographystyle{plain}

\end{document}